\documentclass[11pt]{amsart2000}
\begin{document}
\setlength{\unitlength}{0.01in}
\linethickness{0.01in}
\begin{center}
\begin{picture}(474,66)(0,0)
\multiput(0,66)(1,0){40}{\line(0,-1){24}}
\multiput(43,65)(1,-1){24}{\line(0,-1){40}}
\multiput(1,39)(1,-1){40}{\line(1,0){24}}
\multiput(70,2)(1,1){24}{\line(0,1){40}}
\multiput(72,0)(1,1){24}{\line(1,0){40}}
\multiput(97,66)(1,0){40}{\line(0,-1){40}}
\put(143,66){\makebox(0,0)[tl]{\footnotesize Proceedings of the Ninth Prague Topological Symposium}}
\put(143,50){\makebox(0,0)[tl]{\footnotesize Contributed papers from the symposium held in}}
\put(143,34){\makebox(0,0)[tl]{\footnotesize Prague, Czech Republic, August 19--25, 2001}}
\end{picture}
\end{center}
\vspace{0.25in}
\setcounter{page}{71}
\title{Chainable subcontinua}
\author{Edwin Duda}
\address{University of Miami, Department of Mathematics\\
PO Box 249085\\
Coral Gables, FL 33124-4250}
\email{e.duda@math.miami.edu}
\begin{abstract}
This paper is concerned with conditions under which a metric continuum (a
compact connected metric space) contains a non-degenerate chainable
continuum.
\end{abstract}
\thanks{Edwin Duda,
{\em Chainable subcontinua},
Proceedings of the Ninth Prague Topological Symposium, (Prague, 2001),
pp.~71--73, Topology Atlas, Toronto, 2002}
\maketitle

This paper is concerned with conditions under which a metric continuum (a
compact connected metric space) contains a non-degenerate chainable 
continuum.

By R.H. Bing's theorem eleven \cite{MR13:265a}
if a metric continuum $X$ contains a non-degenerate
subcontinuum $H$ which is hereditarily decomposable, hereditarily
unicoherent, and atriodic, then $H$ is chainable.

The following papers give examples of continua with the property that each
non-degenerate subcontinuum is not chainable. G.T. Whyburn 
\cite{whyburn}.
R.D. Anderson and G. Choquet
\cite{MR21:3819}.
A. Lelek 
\cite{MR26:742}
gives an example of a planar weakly chainable continuum each
non-degenerate subcontinuum of which separates the plane
and thus contains no non-degenerate chainable subcontinuum. 
W.T. Ingram
\cite{MR82k:54056}
gives an example of an hereditarily
indecomposable tree-like continuum such that each non-degenerate
subcontinuum has positive span and hence is not chainable.

C.E. Burgess in 
\cite{MR23:A3551}
shows if a continuum $M$ is almost chainable and $K$ is a proper
subcontinuum of $M$ which contains an endpoint $p$ of $M$, then $K$ is
linearly chainable with $p$ as an end point. 
A continuum $M$ is almost chainable if, for every positive number 
$\varepsilon $, there exists an $\varepsilon$-covering $G$ of $M$ and a
linear chain $C(L_1,L_2,\ldots,L_n)$ of elements of $G$ such that no
$L_i$ $(1\leq i<n)$ intersects an element of $G-C$ and every point of $M$
is within a distance $\varepsilon$ of some element of $C$. 
He also shows if $M$ is almost chainable, then $M$ is not a triod and $M$
is unicoherent and irreducible between some two points. 
Examples show $M$ can contain a triod or a non-unicoherent subcontinuum.

If $X$ and $Y$ are metric continua and if $X$ can be $\varepsilon$-mapped
onto $Y$ for all positive $\varepsilon $ and $Y$ has a non-degenerate
chainable continuum then so does $X$. 
This result suggests considering inverse limit spaces. 
At this stage we refer to a result from the paper of
S. Marde\u{s}i\'{c} and J. Segal
\cite{MR28:1592}
Theorem 1, p.\ 148: ``Every $\pi $-like continuum $X$ is the inverse limit
of an inverse sequence $\{P_i;\pi _{ij}\}$ with bonding maps $\pi _{ij}$
onto and with polyhedra $P_i\in \pi$. 
A continuum is $\pi $-like if it can be $\varepsilon $-mapped onto some
polyhedron in $\pi $ for each positive $\varepsilon $. 
E. Duda and P. Krupski 
\cite{MR93a:54031}
showed that a $k$-junctioned metric continuum, $k$ a non-negative integer,
has at most $k$ points such that any continuum which contains none of the
$k$ points is chainable. A metric continuum is said to be $k$-junctioned
if it is the inverse limit of graphs each of which has at most $k$ branch
points, with surjective bonding maps. A continuum is called finitely
junctioned if it is $k$-junctioned for some non negative integer $k$.

Suppose now $X$ is a tree-like continuum. 
Then for each $\varepsilon >0$ $\;X$ can be mapped onto a tree. 
By a result quoted above $X$ is the inverse limit of a sequence of trees
with surjective bonding maps. 
$X={\displaystyle \lim _{\longleftarrow }}\{T_n,f_{nm} \}$.
Let $f_n:X\rightarrow T_n$ be the standard projection map and let
$$P_n=U\{f^{-1}_n(q)|q \; \mbox{is a branch point}\}.$$ 
Since $T_n$ has at most a finite number of branch points (points of order
$\geq 2$) $P_n$ is closed in $X$. 
If the union of the $P_n$ is not dense in $X$ then $X$ contains a
non-degenerate chainable continuum. 
Actually it is sufficient that $\{P_n\}$ have a subsequence whose union is
not dense in $X$.

Lets now consider a non-degenerate metric continuum in $X$ with span equal
to zero. 
The notion of span was defined by A. Lelek
\cite{MR31:4009}.
In the paper
\cite{MR82c:54031}
he showed continua with span zero are atriodic and tree-like.

There is a series of papers by L.G. Oversteegen and E.D. Tymchatyn which
develop properties of spaces with spans equal to zero or sufficient
conditions that a space have a span equal to zero 
\cite{MR84h:54030, MR86a:54042, MR85j:54051, MR85m:54034}.
Also by L.G. Oversteegen 
\cite{MR91g:54049}.

It is interesting to note that a chainable continuum $X$ can be
$\varepsilon $-mapped onto any fixed dendrite. 
Thus for any tree $T$, by the result of Marde\u{s}i\'{c} and Segal quoted
above, $X$ is the inverse of a sequence of $T$'s.

In the paper 
\cite{MR95c:54056}
P. Minc shows an inverse limit of trees with
simplicial bonding maps having surjective span zero is chainable.

\providecommand{\bysame}{\leavevmode\hbox to3em{\hrulefill}\thinspace}
\providecommand{\MR}{\relax\ifhmode\unskip\space\fi MR }
\providecommand{\MRhref}[2]{%
  \href{http://www.ams.org/mathscinet-getitem?mr=#1}{#2}
}
\providecommand{\href}[2]{#2}

\end{document}